\documentclass[12pt, a4paper]{article}
\usepackage{graphics}
\usepackage{latexsym}
\usepackage{amssymb}
\usepackage{amsmath}

\newcommand{\beq}{\begin{equation}}
\newcommand{\eeq}{\end{equation}}
\newcommand{\beqas}{\begin{eqnarray*}}
\newcommand{\eeqas}{\end{eqnarray*}}
\newcommand{\ep}{\varepsilon}

\newcommand{\spa}{{\mathbb R}^m}

\newcommand{\cb}{{\cal B}}
\newcommand{\cc}{{\cal C}}

\newcommand{\be}{\begin{equation}}
\newcommand{\ee}{\end{equation}}
\newcommand{\bea}{\begin{eqnarray}}
\newcommand{\eea}{\end{eqnarray}}
\newcommand{\beas}{\begin{eqnarray*}}
\newcommand{\eeas}{\end{eqnarray*}}

\def\IP{\mathbb{P}}
\def\IR{\mathbb{R}}

\def\IT{\mathbb{T}}
\def\IE{\mathbb{E}}

\def\dw{\dot{W}}

\def\ca{{\cal A}^W}
\def\cab{{\cal A}^{W,b}}

\def\cg{{\cal G}}
\def\cn{{\cal N}}

\newtheorem{theorem}{Theorem}[section]
\newtheorem{lemma}[theorem]{Lemma}
\newtheorem{assumption}[theorem]{Assumption}
\newtheorem{proposition}[theorem]{Proposition}
\newtheorem{definition}[theorem]{Definition}
\newtheorem{corollary}[theorem]{Corolarry}

\numberwithin{equation}{section}
\begin{document}

\title{\bf Bayesian inverse problems for Burgers and Hamilton-Jacobi equations with white noise forcing}
\author{Viet Ha Hoang\\Division of Mathematical Sciences,\\
School of Physical and Mathematical Sciences,\\
Nanyang Technological University, Singapore 637371}

\date{}
\maketitle

\begin{abstract}
The paper formulates Bayesian inverse problems for inference in a topological measure space given noisy observations. Conditions for the validity of the Bayes formula and the well-posedness of the posterior measure are studied. The abstract theory is then applied to Burgers and Hamilton-Jacobi equations on a semi-infinite time interval with forcing functions which are white noise in time. Inference is made on the white noise forcing, assuming the Wiener measure as the prior. 
\end{abstract}
\section{Introduction}
Bayesian inverse problems are attracting considerable attention due to their importance in many applications (\cite{Stuart}). Given a functional $\cg:X\to\spa$ where $X$ is a probability space, the observation $y\in\IR^m$ of $\cg$ is subject to an unbiased noise $\sigma$:
\be
y=\cg(x)+\sigma.
\label{eq:basic}
\ee
The inverse problem determines $x$ given the noisy observation $y$. 

Cotter et al. \cite{CDRS} formulate a rigorous mathematical framework for Banach spaces $X$, assuming a Gaussian prior probability measure. They prove the Bayes formula for this infinite dimensional setting, and determine the Radon-Nikodym derivative of the posterior measure with respect to the prior. When $\cg(x)$ grows polynomially, using Fernique theorem, they show that the posterior measure continuously depends on $y$; the problem is thus well-posed. 

Cotter et al. \cite{CDRS} then apply the framework to Bayesian inverse problems for Navier Stokes equation for data assimilation in fluid mechanics. They consider  Eulerian and Lagrangian data assimilation where the fluid velocity and traces of passive floats are observed respectively. They make inference on the initial velocity, and in the case of model errors, on the forcing. The random forcing function is  an Ornstein–-Uhlenbeck process which defines a Gaussian probability measure on a Hilbert function space. Computationally, the posterior measure is sampled via the Markov Chain-Monte Carlo method in Cotter et al. \cite{Cotteretalfluid}.  

In this paper, we are interested in Bayesian inverse problems for partial differential equations where forcing is white noise in time. We recover the white noise forcing on a semi-infinite time interval $(-\infty, t]$, given noisy observations.  Semi-infinite time intervals are interested when the long time behaviour of the solutions is desired. We therefore depart from the Banach space framework of Cotter et al. \cite{CDRS} and consider the measure space setting instead. Our first aim is to formulate Bayesian inverse problems in measure spaces, for which we establish the validity of the Bayes formula and the well-posedness of the posterior measure. Our second aim is to apply the theory to inference on partial differential equations with a forcing function which is white noise in time, in particular, we study the randomly forced Burgers and Hamilton-Jacobi equations.

Burgers equation is a well known model for turbulence (see, e.g., \cite{BecKhanin}). It is simpler than Navier-Stokes equation due to the absence of the pressure but still keeps the essential nonlinear features. It also arises in many problems in non-equilibrium statistical mechanics. Burgers equation and equations of the Burgers type are employed as models for studying data assimilation in Bardos and Pironneau \cite{BardosPironneau}, Lundvall et al. \cite{Lundvalletal}, Apte et al. \cite{Apteetal} and Freitag et al. \cite{Freitagetal}. We consider the case where spatially, the forcing is the gradient of a potential. Assuming that the velocity is also a gradient, the velocity potential satisfies a Hamilton-Jacobi equation with random forcing. This is the Kardar-Parisi-Zhang (KPZ) equation that describes lateral growth of an interface. Burgers equation with gradient forcing and Hamilton Jacobi equation are thus closely related. In the dissipative case, inverse problems to determine the dissipative coefficient for the KPZ equation was studied in Lam and Sander \cite{LS}. We will only concentrate on inviscid problems. Inviscid Burgers equation possesses many shocks where the solution is not continuous.

Burgers equation with white noise forcing in the spatially periodic case is studied by E et al. \cite{EKMS} in one dimension and by Iturriaga and Khanin \cite{IK} in multidimensions. They prove that almost surely, there is a unique solution that exists for all time.  They establish many properties for shocks and for the dynamics of the Lagrange minimizers of the Lax-Oleinik functional. Much less is understood in the non-compact case. Hoang and Khanin \cite{HoangKhanin} study the special case where the forcing potential has a large maximum and a small minimum. They prove that there is a solution that exists for all time. It is the limit of finite time solutions with the zero initial condition. 

For  Hamilton-Jacobi equation \eqref{eq:HJ}, we assume that the solution (i.e. the velocity potential $\phi$) is observed at fixed spatial points at fixed moments. This assumption is reasonable as $\phi$ is continuous and is well defined  everywhere (though it is determined within an additive constant). Solutions for Burgers equation \eqref{eq:Burgers}  are not defined at shocks, but for a given time, they are in $L^1_{\rm loc}$, so observations in the form of a continuous functional in $L^1_{\rm loc}$ are well-defined. 

In section 2, we establish the Bayes formula for the posterior measures  of Bayesian inverse problems in measure spaces, and propose conditions under which the problems are well-posed. Polynomial growth and Fernique theorem are no longer used as we do not assume a Gaussian prior. The results are generalizations of those in the Banach space setting by Cotter et al. in \cite{CDRS}. Section 3 introduces the randomly forced Burgers and Hamilton-Jacobi equations and  formulates the Bayesian inverse problems. The spatially periodic case is considered in section 4, using the results for random Burgers and Hamilton-Jacobi equations by E et al. \cite{EKMS} and Iturriaga and Khanin \cite{IK}. We establish the validity of the Bayes formula and the posterior measure's local Lipschitzness with respect to the Hellinger metric; the main results are recorded in Theorems \ref{thm:HJp} and \ref{thm:Burgersp}. The non-periodic case is studied in section 5, using the results by Hoang and Khanin \cite{HoangKhanin} where we restrict our consideration to  forcing potentials  that satisfy Assumption \ref{assum:np}, with a large maximum and a small minimum. Under this assumption, for both Hamilton-Jacobi and Burgers equations, the Bayes formula is valid and the posterior measure is locally Lipschitz with respect to the Hellinger metric. The main results are presented in Theorem \ref{thm:HJnp} and \ref{thm:Burgersnp}.  Section 6 proves some technical results that are used in the previous two sections. 

Throughout the paper, we define $\langle \cdot,\cdot\rangle_\Sigma=\langle\Sigma^{-1/2}\cdot,\Sigma^{-1/2}\cdot\rangle$ for a positive definite matrix $\Sigma$; the corresponding norm $|\Sigma^{-1/2}\cdot|$ is denoted as $|\cdot|_\Sigma$. We denote by $c$ various positive constants that do not depend on the white noise; the value of $c$ can differ from one appearance to the next.

\section{Bayesian inverse problems in measure spaces}\label{sec:Bayesian}
We study Bayesian inverse problems in measure spaces, in particular, the conditions for the validity of the Bayes formula, and the well-posedness of the problems. In \cite{CDRS}, Cotter et al. require the functional $\cg$ in \eqref{eq:basic} to be measurable with respect to the prior measure for the validity of the Bayes formula. Restricting their consideration to Banach spaces, this requirement holds when $\cg$ is continuous.  However, this is true for any topological spaces $X$. Cotter et al. imposed a polynomial growth on $\cg$ to prove the well-posedness. Using Fernique theorem, they show that when the prior measure is Gaussian, with respect to the Hellinger metric, the posterior measure is locally Lipschitz with respect to $y$. We will generalize this result to measure spaces, and introduce a general condition for the local Lipschitzness of the posterior measure, without imposing a Gaussian prior.  

\subsection{Problem formulation}
For a topological space $X$ with a $\sigma$-algebra ${\cal B}(X)$, we consider a functional  ${\cal G}:X\to {\mathbb R}^m$. An observation $y$ of ${\cal G}$ is subject to a multivariate Gaussian noise $\sigma$, i.e
\be
y=\cg(x)+\sigma,
\label{eq:prob}
\ee
where $\sigma\sim\cn(0,\Sigma)$. 

Let $\mu_0$ be the prior probability measure on $(X,{\cal B}(X))$. Our purpose is to determine the conditional probability $\IP(x|y)$.  Let 
\be
\Phi(x,y)={1\over 2}|y-\cg(x)|^2_{\Sigma}.
\label{eq:density}
\ee
We have the following theorem.
\begin{theorem}(Cotter et al. \cite{CDRS} Theorem 2.1)
If $\cg:X\to\spa$ is $\mu_0$-measurable then the posterior measure $\mu^y(dx)=\IP(dx|y)$ is absolutely continuous with respect to the prior measure $\mu_0(dx)$ and the Radon-Nikodym derivative satisfies
\beq
{d\mu^y\over d\mu_0}(x)\propto\exp(-\Phi(x;y)).\label{eq:RN}
\eeq
\end{theorem}
Although Cotter et al. assume that $X$ is a Banach space, their proof  holds for any topological measure spaces $X$. 

Following \cite{CDRS} corollary 2.2, we have the following corollary.
\begin{corollary}\label{corollary:2.1} If $\cg:X\to\spa$ is continuous, then the measure $\mu^y(dx)=\IP(dx|y)$ is absolutely continuous with respect to $\mu_0(x)$ and the Radon-Nikodym derivative satisfies (\ref{eq:RN}).
\end{corollary}
{\it Proof}\ \ If $\cg:X\to\spa$ is continuous, then it is measurable (\cite{Kallenberg}, Lemma 1.5). The conclusion holds. \hfill$\Box$

\subsection{Well-posedness of the problem}
For Banach spaces $X$, when the prior measure $\mu_0$ is Gaussian, the well-posedness of the inverse problem is established by Cotter et al. \cite{CDRS} by assuming that the function $\Phi(x;y)$ has a polynomial growth in $\|x\|_X$ fixing $y$, and is Lipschitz in $y$ with the Lipschitz coefficient being bounded by a polynomial of $\|x\|_X$. The proof uses Fernique theorem. For a general measure space $X$, we make the following assumption.
\begin{assumption} \label{assumption}
The function $\Phi:X\times \spa\to {\mathbb R}$ satisfies:
\begin{itemize}

\item[(i)] For each $r>0$, there is a constant $M(r)>0$ and a set $X(r)\subset X$ of positive $\mu_0$ measure such that for all $x\in X(r)$ and for all $y$ such that $|y|\le r$
\[
0\le\Phi(x;y)\le M(r).
\]
\item[(ii)] There is a function $G:{\mathbb R}\times X\to {\mathbb R}$ such that for each $r>0$, $G(r,\cdot)\in L^2(X,d\mu_0)$ and if $|y|\le r$ and $|y'|\le r$ then
\[
|\Phi(x;y)-\Phi(x;y')|\le G(r,x)|y-y'|.
\]
\end{itemize}
\end{assumption}

To study the well-posedness of the Bayesian inverse problems (\ref{eq:prob}), following Cotter et al. \cite{CDRS}, we employ the Hellinger metric 
\[
d_{\rm{Hell}}(\mu,\mu')=\sqrt{{1\over 2}\int_X\Bigl(\sqrt{{d\mu\over d\mu_0}}-\sqrt{{d\mu'\over d\mu_0}}\Bigr)^2d\mu_0},
\]
where $\mu$ and $\mu'$ are measures on $(X,\cb(X))$.
The probability measure $\mu^y$ is determined as
\[
{d\mu^y\over d\mu_0}={1\over Z(y)}\exp(-\Phi(x,y)),
\]
where the normalization constant is
\be
Z(y)=\int_X\exp(-\Phi(x,y))d\mu_0(x).
\label{eq:zy}
\ee
With respect to the Hellinger metric, the measure $\mu^y$ depends continuously on  $y$ as shown in the following.
\begin{theorem}\label{thm:wellposedness} Under Assumption \ref{assumption}, the measure $\mu^y$ is locally Lipschitz in the data $y$ with respect to the Hellinger metric: for each positive constant $r$ there is a positive constant $c(r)$ such that if $|y|\le r$ and $|y'|\le r$, then
\[
d_{\rm{Hell}}(\mu^y,\mu^{y'})\le c(r)|y-y'|.
\]
\end{theorem}
 {\it Proof} \ \ 
First we show that for each $r>0$, there is a positive constant $c(r)$ such that $Z(y)\ge c(r)$ when $|y|\le r$. It is obvious from (\ref{eq:zy}) and from Assumption \ref{assumption}(i) that when $|y|\le r$:
\be
Z(y)\ge \mu_0(X(r))\exp(-M(r)).
\label{eq:boundforZ}
\ee
Using the inequality $|\exp(-a)-\exp(-b)|\le |a-b|$ for $a>0$ and $b>0$, we have
\[
|Z(y)-Z(y')|\le \int_X|\Phi(x;y)-\Phi(x;y')|d\mu_0(x).
\]
From Assumption \ref{assumption}(ii), when $|y|\le r$ and $|y'|\le r$:
\[
|\Phi(x;y)-\Phi(x;y')|\le G(r,x)|y-y'|.
\]
As $G(r,x)$ is $\mu_0$-integrable, 
\be
|Z(y)-Z(y')|\le c(r)|y-y'|.
\label{eq:LipZ}
\ee
The Hellinger distance satisfies
\beqas
2d_{\rm{Hell}}(\mu^y,\mu^{y'})^2&=&\int_X\Bigl(Z(y)^{-1/2}\exp(-{1\over 2}\Phi(x,y))-\\
&&{\hskip 2cm}Z(y')^{-1/2}\exp(-{1\over 2}\Phi(x,y'))\Bigr)^2d\mu_0(x)\\
&\le& I_1+I_2,
\eeqas
where 
\[
I_1={2\over Z(y)}\int_X\Bigl(\exp(-{1\over 2}\Phi(x,y))-\exp(-{1\over 2}\Phi(x,y'))\Bigr)^2d\mu_0(x),
\]
and
\[
I_2=2|Z(y)^{-1/2}-Z(y')^{-1/2}|^2\int_X\exp(-\Phi(x,y'))d\mu_0(x).
\]
Using again the inequality $|\exp(-a)-\exp(-b)|\le|a-b|$, we deduce from \eqref{eq:boundforZ}
\beqas
I_1&\le& c(r)\int_X|\Phi(x,y)-\Phi(x,y')|^2d\mu_0(x)\\
&\le&c(r)\int_X(G(r,x))^2d\mu_0(x)|y-y'|^2\le C(r)|y-y'|^2.
\eeqas 
Furthermore,
\[
|Z(y)^{-1/2}-Z(y')^{-1/2}|^2={|Z(y)-Z(y')|^2\over (Z(y)^{1/2}+Z(y')^{1/2})Z(y)Z(y')}.
\]
From \eqref{eq:boundforZ} and \eqref{eq:LipZ},
\[
|Z(y)^{-1/2}-Z(y')^{-1/2}|^2
\le c(r)|Z(y)-Z(y')|^2\le c(r)|y-y'|^2.
\]
The conclusion then follows.\hfill$\Box$

\section{Bayesian inverse problems for equations with white noise forcing}
\subsection{Stochastic Burgers and Hamilton-Jacobi equations}
We consider the inviscid Burgers equation 
\beq
{\partial u\over\partial t}+(u\cdot\nabla)u=f^W(x,t)
\label{eq:Burgers}
\eeq
where $u(x,t)\in\IR^d$ is the velocity of a fluid particle in $\IR^d$.
The forcing function
$f^{W}(x,t)\in\IR^d$ depends on $t$ via a one dimensional white noise $\dot W$ and is of the form
\[
f^{W}(x,t)=f(x)\dot{W}(t).
\]
We study the case where $f(x)$ is a gradient, i.e. there is a potential $F(x)$ such that
\[
f(x)=-\nabla F(x).
\]
Throughout this paper, $F(x)$ is  three time differentiable with bounded derivatives up to the third order. We
consider the case where $u(x,t)$ is also a gradient, i.e. $u(x,t)=\nabla\phi(x,t)$. The velocity potential $\phi$ satisfies the Hamilton-Jacobi equation
\beq
{\partial\phi(x,t)\over\partial t}+{1\over 2}|\nabla\phi(x,t)|^2+F^W(x,t)=0,
\label{eq:HJ}
\eeq
where the forcing  potential $F^W(x,t)$ is
\[
F^W(x,t)=F(x)\dw(t).
\]
 
The viscosity solution $\phi$ of (\ref{eq:HJ}) that corresponds to the solution $u$ of the Burgers equation \eqref{eq:Burgers} is determined by the Lax-Oleinik variational principle. With the initial condition $\phi_0$ at time 0, $\phi(x,t)$ is given by
\beq
\phi(x,t)=\inf\biggl\{\phi_0(\gamma(0))+\int_0^t\Bigl({1\over 2}|{\dot\gamma}|^2-F^W(\gamma(\tau),\tau)\Bigr)d\tau\biggr\},
\label{eq:LO}
\eeq
where the infimum is taken among all the absolutely continuous curves $\gamma:[0,t]\to\spa$ such that $\gamma(t)=x$.

When the forcing function is spatially periodic,  problem (\ref{eq:HJ}) possesses a unique solution (up to an additive constant) that exists for all time $t$ (see \cite{EKMS,IK}). The non-compact setting is far more complicated. In (\cite{HoangKhanin}), it is proved that when $F(x)$ has a large maximum, and a small minimum, (\ref{eq:HJ}) possesses a solution that exists for all time. These imply the existence of a solution to the Burgers equation that exists for all time.

We employ the framework for measure spaces developed in section \ref{sec:Bayesian} to study Bayesian inverse problems for Burgers equation \eqref{eq:Burgers} and Hamilton-Jacobi equation \eqref{eq:HJ} with white noise forcing.

\subsection{Bayesian inverse problems}

We study the Bayesian inverse problems that make inference on  the white noise forcing given the observations at a finite set of times. 

For Hamilton-Jacobi equation \eqref{eq:HJ},
 fixing $m$ spatial points $x_i$ and $m$ times $t_i$, $i=1,\ldots,m$, the function $\phi$ is observed at $(x_i,t_i)$. The observations are subject  to an unbiased Gaussian noise. They are
\[
z_{i}=\phi(x_i,t_i)+\sigma'_{i},
\]
where $\sigma_{i}'$ are independent Gaussian random variables.

Since $\phi$ is determined within an additive constant, we assume that it is measured at a further point $(x_0,t_0)$ where, for simplicity only,  $t_0<t_i$ for all $i=1,\ldots,m$. 

Let
\[
y_{i}=z_{i}-z_{0}=\phi(x_i,t_i)-\phi(x_0,t_0)+\sigma'_{i}-\sigma'_{0}.
\]
Let $\sigma_{i}=\sigma'_{i}-\sigma'_{0}$. We assume that $\sigma=\{\sigma_{i}:\ i=1,\ldots,m\}\in\IR^m$ follows a Gaussian distribution with zero mean and covariance $\Sigma$.  The vector 
\be
{\cal G}_{HJ}(W)=\{\phi(x_i,t_i)-\phi(x_0,t_0):\ i=1,\ldots,m\}\in\spa
\label{eq:GW}
\ee
 is uniquely determined by the Wiener process $W$. We recover $W$ given 
\[
y=\cg_{HJ}(W)+\sigma.
\]
The solution $u(x,t)$ of  Burgers equation \eqref{eq:Burgers} is not defined at shocks. However, in the cases studied here, $u(\cdot,t)\in L^1_{\rm{loc}}(\spa)$ for all time $t$. We assume that noisy measurements 
\[
y_i=l_i(u(\cdot,t_i))+\sigma_i,\ i=1,\ldots,m,
\]
are made for $l_i(u(\cdot,t_i))$ ($i=1,\ldots,m$) where $l_i:L^1_{\rm{loc}}(\spa)\to \IR$ are continuous and bounded functionals. The noise $\sigma=(\sigma_1,\ldots,\sigma_m)\sim {\cal N}(0,\Sigma)$. 
Letting 
\be
{\cal G}_B(W)=(l_1(u(\cdot,t_1)),\ldots,l_m(u(\cdot,t_m)))\in\spa,
\label{eq:GWB}
\ee
given
\[
y=\cg_B(W)+\sigma,
\]
we determine $W$.

Let $t_{\max}=\max\{t_1,\ldots,t_m\}$.
Let $\cal X$ be the metric space of continuous functions $W$ on $(-\infty,t_{\max}]$ with $W(t_{\max})=0$. The space is equipped with the metric:
\beq
D(W,W')=\sum_{n=1}^\infty{1\over 2^n}{\sup_{-n\le t\le t_{\max}}|W(t)-W'(t)|\over 1+\sup_{-n\le t\le t_{\max}}|W(t)-W'(t)|}
\label{eq:metric}
\eeq
(see Stroock and Varadhan \cite{SV}).
Let $\mu_0$ be the Wiener measure on $\cal X$. Let $\mu^y$ be the conditional probability defined on $\cal X$ given $y$. For $\cg_{HJ}$ in \eqref{eq:GW}, we define the function $\Phi_{HJ}:{\cal X}\times \IR^m\to \IR$ by
\be
\Phi_{HJ}(W;y)={1\over 2}|y-{\cal G}_{HJ}(W)|^2_{\Sigma},
\label{eq:PhiHJ}
\ee
and for $\cg_B$ in \eqref{eq:GWB}, we define
\be
\Phi_{B}(W;y)={1\over 2}|y-{\cal G}_{B}(W)|^2_{\Sigma}.
\label{eq:PhiB}
\ee

We will prove that the Radon-Nikodym derivative of $\mu^y$ satisfies
\[
{d\mu^y\over d\mu_0}\propto\exp(-\Phi_{HJ}(W;y)),
\]
and 
\[
{d\mu^y\over d\mu_0}\propto\exp(-\Phi_{B}(W;y))
\]
respectively, 
and that the posterior measure $\mu^y$ continuously depends on $y\in\spa$.

\section{Periodic problems}
We first consider the case where the problems are set in the torus $\IT^d$. The forcing potential $F^{W}(x,t)$ is spatially periodic. Randomly forced Burgers and Hamilton-Jacobi equations in a torus are  studied thoroughly by E et al. in \cite{EKMS} for one dimension and by Iturriaga and Khanin in \cite{IK} for multidimensions. We first review their results. 
\subsection{Existence and uniqueness for periodic problems}
The mean value
\[
b=\int_{\IT^d} u(x,t)dx
\]
is unchanged for all $t$. We denote the solution $\phi$ of the Hamilton-Jacobi equation \eqref{eq:HJ} by $\phi^W_b$. It is of the form
\beq
\phi^W_b(x,t)=b\cdot x+\psi^W_b(x,t),
\label{psi}
\eeq
where $\psi^W_b$ is spatially periodic. Let $\cc(s,t;\IT^d)$ be the space of absolutely continuous functions from $(s,t)$ to $\IT^d$. 
For each vector $b\in{\mathbb R}^d$, we define the operator $\cab_{s,t}:\cc(s,t;\IT^d)\to\IR$ as
\be
\cab_{s,t}(\gamma)=\int_s^t\Bigl({1\over 2}|{\dot\gamma}(\tau)-b|^2-F^W(\gamma(\tau),\tau)-{b^2\over 2}\Bigr)d\tau.
\label{eq:AWb}
\ee 
The variational formulation for the function $\psi^W_b(x,t)$ in (\ref{psi}) is
\be
\psi^W_b(x,t)=\inf_{\gamma\in\cc(s,t;\IT^d)}\Bigl(\psi(\gamma(s),s)+\cab_{s,t}(\gamma)\Bigr),
\label{eq:Lax}
\ee
where $s<t$. 
This is written in terms of the Lax operator as
\[
{\cal K}_{s,t}^{W,b}\psi^W_b(\cdot,s)=\psi^W_b(\cdot,t).
\]
Using integration by parts, we get
\begin{eqnarray}
\cab_{s,t}(\gamma)=\int_s^t\biggl({1\over 2}|{\dot\gamma}(\tau)-b|^2+\nabla F(\gamma(\tau))\cdot{\dot\gamma}(\tau)(W(\tau)-W(t))-{b^2\over 2}\biggr)d\tau\nonumber\\
-F(\gamma(s))(W(t)-W(s)).\quad
\label{eq:AWb1}
\end{eqnarray}
We study the functional $\cab_{s,t}$ via minimizers which are defined as follows (\cite{IK}).
\begin{definition}\label{def:minimizer}
A curve $\gamma:[s,t]\to \IR^d$ such that $\gamma(s)=x$ and $\gamma(t)=y$ is called a minimizer over $[s,t]$ if it minimizes the action $\cab_{s,t}$ among all the absolutely continuous curves with end points at $x$ and $y$ at times $s$ and $t$ respectively. 

For a function $\psi$, a curve $\gamma:[s,t]\to \IR^d$ with $\gamma(t)=x$ is called a $\psi$ minimizer over $[s,t]$ if it minimizes the action $\psi(\gamma(s))+\cab_{s,t}(\gamma)$ among all the curves with end point at $x$ at time $t$.

A curve $\gamma:(-\infty,t]\to \IR^d$ with $\gamma(t)=x$ is called a one-sided minimizer if it is a minimizer over all the time intervals $[s,t]$.
\end{definition} 

Iturriaga and Khanin \cite{IK} proved the following theorem :
\begin{theorem}\label{thm:periodic}(\cite{IK}, Theorem 1)
(i) For almost all $W$, all $b\in \IR^d$ and  $T\in\IR$ there exists a unique (up to an additive constant) function $\phi_b^W(x,t)$, $x\in\spa$, $t\in(-\infty,T]$ such that
\[
\phi_b^W(x,t)=b\cdot x+\psi_b^W(x,t),
\]
where $\psi_b^W(\cdot,t)$ is a $\IT^d$ periodic function and for all $-\infty<s<t<T$,
\[
\psi_b^W(\cdot,t)={\cal K}_{s,t}^{W,b}\psi_b^W(\cdot,s).
\]
(ii) The function $\psi_b^W(\cdot,t)$ is Lipschitz in $\IT^d$. If $x$ is a point of differentiability of $\psi_b^W(\cdot,t)$ then there exists a unique one-sided minimizer $\gamma_{x,t}^{W,b}$ at $(x,t)$ and its velocity is given by the gradient of $\phi$: ${\dot\gamma}^{W,b}_{x,t}(t)=\nabla\phi_b^W(x,t)=b+\nabla\psi_b^W(x,t)$. Further, any onesided minimizer is a $\psi^W_b$ minimizer on finite time intervals. 

(iii) For almost all $W$ and all $b\in\IR^d$:
\[
\lim_{s\to-\infty}\sup_{\eta\in C(\IT^d)}\min_{C\in\IR}\max_{x\in\IT^d}|{\cal K}^{W,b}_{s,t}\eta(x)-\psi_b^W(x,t)-C|=0.
\]
(iv) The unique solution $u^W_b(x,t)$ that exists for all time of the Burgers equation \eqref{eq:Burgers} is determined by $u^W_b(\cdot,t)=\nabla\phi^W_b(\cdot,t)$.
\end{theorem}

We will employ these results to study the Bayesian inverse problems formulated in section \ref{sec:Bayesian}.

\subsection{Bayesian inverse problem for spatially periodic Hamilton-Jacobi equation \eqref{eq:HJ}}
We first study the Bayesian inverse problem for Hamilton-Jacobi equation \eqref{eq:HJ} in the torus $\IT^d$.
The observation  ${\cal G}_{HJ}(W)$ in \eqref{eq:GW} becomes 
\[
{\cal G}_{HJ}(W)=\{\psi_b^W(x_i,t_i)-\psi_b^W(x_0,t_0):\ i=1,\ldots,m\}.
\]
To establish the validity of the Bayesian formula \eqref{eq:RN}, we first show that $\cg_{HJ}(W)$ as a map from $\cal X$ to $\IR^m$ is continuous with respect to the metric $D$ in \eqref{eq:metric}.

The following Proposition holds.
\begin{proposition}\label{prop:concghjp} The map ${\cal G}_{HJ}(W):{\cal X}\to {\mathbb R}^m$ is continuous with respect to the metric \eqref{eq:metric}.
\end{proposition}
{\it Proof}\ \ 
Let $W_k$ converge to $W$ in the metric \eqref{eq:metric}.
There are constants $C_k$ which do not depend on $i$  such that  
\beq
\lim_{k\to\infty}|\phi_b^{W_k}(x_i,t_i)- \phi_b^{W}(x_i,t_i)-C_k|\to 0
\label{eq:epperiodic}
\eeq
for all $i=0,\ldots,m$. 
The proof of \eqref{eq:epperiodic} is technical; we present it in section \ref{sec:6.1}.
From this we deduce
\[
\lim_{k\to\infty}|(\phi_b^{W_k}(x_i,t_i)-\phi_b^{W_k}(x_0,t_0))-(\phi_b^{W}(x_i,t_i)-\phi_b^{W}(x_0,t_0))|=0,
\]
so
\[
\lim_{k\to\infty}|{\cal G}_{HJ}(W_k)-{\cal G}_{HJ}(W)|=0.
\]
The conclusion then follows. \hfill$\Box$

To establish the continuous dependence of $\mu^y$ on $y$, we verify Assumption \ref{assumption}. First we prove a bound for $\cg_{HJ}(W)$. 
\begin{lemma}
There is a positive constant $c$ which only depends on the potential $F$ such that
\[
|\cg_{HJ}(W)|\le c\Big(1+\sum_{i=1}^m\max_{t_0-1\le \tau\le t_i}|W(\tau)-W(t_i)|^2\Big).
\]
\end{lemma}
{\it Proof}\ \ 
Let $\gamma$ be a $\psi^W_b(\cdot,t_0-1)$ minimizer starting at $(x_i,t_i)$. Then 
\[
\psi^W_b(x_i,t_i)=\psi_b^W(\gamma(t_0-1),t_0-1)+\cab_{t_0-1,t_i}({\gamma}).
\]
From \eqref{eq:AWb1}, there is a constant $c$ which depends only on $b$ and $F$ such that 
\beqas
\cab_{t_0-1,t_i}(\gamma)&\ge& {1\over 2(t_i-t_0+1)}\Bigl(\int_{t_0-1}^{t_i}|\dot\gamma(\tau)|d\tau\Bigr)^2-\\
&&c(1+\max_{t_0-1\le\tau\le t_i}|W(\tau)-W(t_i)|)\int_{t_0-1}^{t_i}|\dot\gamma(\tau)|d\tau\\
&&-c(1+|W(t_i)-W(t_0-1)|.
\eeqas
The right hand side is a quadratic form of $\int_{t_0-1}^{t_i}|\dot\gamma(\tau)|d\tau$ so
\[
\cab_{t_0-1,t_i}(\gamma)\ge -c(1+\max_{t_0-1\le\tau\le t_i}|W(\tau)-W(t_i)|^2).
\]
Therefore
\[
\psi^W_b(x_i,t_i)\ge\psi_b^W(\gamma(t_0-1),t_0-1)-c(1+\max_{t_0-1\le\tau\le t_i}|W(\tau)-W(t_i)|^2).
\]
Let $\gamma'$ be the linear curve connecting $(x_0,t_0)$ and $(\gamma(t_0-1),t_0-1)$. We have
\beqas
\psi^W_b(x_0,t_0)&\le& \psi^W_b(\gamma(t_0-1),t_0-1)+\cab_{t_0-1,t_0}(\gamma')\\
&&\le \psi_b^W(\gamma(t_0-1),t_0-1)+c(1+\max_{t_0-1\le\tau\le t_0}|W(\tau)-W(t_0)|).
\eeqas
Therefore 
\be
\psi_b^W(x_i,t_i)-\psi_b^W(x_0,t_0)\ge -c(1+\max_{t_0-1\le\tau\le t_i}|W(\tau)-W(t_i)|^2).
\label{eq:eq1}
\ee 
Let ${\bar\gamma}$ be the linear curve connecting $(x_0,t_0)$ and $(x_i,t_i)$. We then find that
\[
\psi_b^W(x_i,t_i)\le \psi_b^W(x_0,t_0)+\cab_{t_0,t_i}({\bar\gamma}).
\]
Thus 
\be
\psi_b^W(x_i,t_i)\le\psi_b^W(x_0,t_0)+c(1+\max_{t_0\le\tau\le t_i}|W(\tau)-W(t_i)|).
\label{eq:eq2}
\ee
From \eqref{eq:eq1} and \eqref{eq:eq2}, we conclude that there exists a positive constant $c$ such that 
\be
|\psi_b^W(x_i,t_i)-\psi_b^W(x_0,t_0)|\le c(1+\max_{t_0-1\le\tau\le t_i}|W(\tau)-W(t_i)|^2).
\ee
Therefore
\be
|\cg_{HJ}(W)|\le c\Big(1+\sum_{i=1}^m\max_{t_0-1\le\tau\le t_i}|W(\tau)-W(t_i)|^2\Big).
\label{eq:boundp}
\ee
\hfill$\Box$

\begin{proposition}\label{prop:assumHJp} Assumption \ref{assumption} holds.
\end{proposition}
{\it Proof}\ \ Fixing a constant $M$, from \eqref{eq:boundp} there is a set of positive Wiener measure so that 
\[
|\cg_{HJ}(W)|\le M.
\]
 Therefore, there is a constant $M(r)$ and a set of positive measure such that when $|y|\le r$
\[
|\Phi_{HJ}(W;y)|={1\over 2}|y-\cg_{HJ}(W)|_{\Sigma}^2\le {1\over 2}\|\Sigma^{-1/2}\|_{\spa,\spa}^2(|y|+M)^2\le M(r).
\]
Assumption \ref{assumption}(i) holds.

Now we prove Assumption \ref{assumption}(ii). We have
\begin{eqnarray}
|\Phi_{HJ}(W;y)-\Phi_{HJ}(W;y')|= {1\over 2}\Big|\langle\Sigma^{-1/2}(y-\cg_{HJ}(W)),\Sigma^{-1/2}(y-y')\rangle+\nonumber\\
\langle\Sigma^{-1/2}(y-y'),\Sigma^{-1/2}(y'-\cg_{HJ}(W))\rangle\Big|\nonumber\\
\le {1\over 2}\|\Sigma^{-1/2}\|_{\spa,\spa}^2(|y|+|y'|+2|\cg_{HJ}(W)|)|y-y'|.\qquad\label{eq:Lipy}
\end{eqnarray}
Letting
\[
G(r,W)=\|\Sigma^{-1/2}\|_{\spa,\spa}^2\bigg[r+c\Big(1+\sum_{i=1}^m\max_{t_0-1\le\tau\le t_i}|W(\tau)-W(t_i)|^2\Big)\bigg].
\]
Obviously, $G(r,\cdot)\in L^2({\cal X},\mu_0)$. 
\hfill$\Box$

From Propositions \ref{prop:concghjp}, \ref{prop:assumHJp}, Corollary \ref{corollary:2.1} and Theorem \ref{thm:wellposedness} we have:
\begin{theorem}\label{thm:HJp} For spatially periodic Hamilton-Jacobi equation, the measure $\mu^y(dW)=\IP(dW|y)$ is absolutely continuous with respect to the Wiener measure $\mu_0(dW)$; the Radon-Nikodym derivative satisfies
\[
{d\mu^y\over d\mu_0}\propto \exp(-\Phi_{HJ}(W;y)).
\]
When $|y|\le r$ and $|y'|\le r$, there is a constant $c(r)$ such that
\[
d_{\rm Hell}(\mu^y,\mu^{y'})\le c(r)|y-y'|.
\]
\end{theorem}

\subsection{Bayesian inverse problem for spatially periodic Burgers equation}
We consider the Bayesian inverse problem for Burgers equation \eqref{eq:Burgers} in the torus $\IT^d$. We first establish the continuity of 
\[
\cg_B(W)=(l_1(u^W_b(\cdot,t_1)),\ldots,l_m(u^W_b(\cdot,t_m)))
\]
as a map from $\cal X$ to $\IR^m$ to prove the validity of the Bayesian formula \eqref{eq:RN}. 

First we prove the uniform boundedness of $|\dot\gamma(t)|$ for any minimizers $\gamma$ on $[t-1,t]$. The result is proved in \cite{EKMS} and \cite{IK}, using Gronvall's inequality. We present another proof without using this inequality, and establish explicitly the dependency of the bound on $W$.
\begin{lemma}\label{lem:dotgammat}
For any minimizer $\gamma$ of $\cab_{t-1,t}$,
\[
|\dot\gamma(t)|\le c(1+\max_{t-1\le\tau\le t}|W(\tau)-W(t)|^2).
\]
\end{lemma}
{\it Proof}\ \
Let $\bar\gamma$ be the linear curve connecting $(\gamma(t),t)$ and $(\gamma(t-1),t-1)$. As $|\gamma(t)-\gamma(t-1)|\le d^{1/2}$, 
\[
\cab_{t-1,t}(\bar\gamma)\le c(1+\max_{t-1\le\tau\le t}|W(\tau)-W(t)|).
\]
 We have
\beqas
\cab_{t-1,t}(\gamma)\ge {1\over 2}\Big(\int_{t-1}^t|\dot\gamma(\tau)|d\tau\Big)^2-c(1+\max_{t-1\le\tau\le t}|W(\tau)-W(t)|)\int_{t-1}^t|\dot\gamma(\tau)|d\tau-\\
c(1+|W(t-1)-W(t)|).
\eeqas
Because $\cab_{t-1,t}(\gamma)\le \cab_{t-1,t}(\bar\gamma)$, solving the quadratic equation, we get
\beq
\int_{t-1}^t|\dot\gamma(\tau)|d\tau\le c(1+\max_{t-1\le\tau\le t}|W(\tau)-W(t)|).
\label{eq:dotgammat}
\ee
The minimizer $\gamma$ satisfies (see \cite{EKMS} and \cite{IK})
\begin{eqnarray}
\dot\gamma(t)=\dot\gamma(s)+\int_s^t\nabla^2F(\gamma(\tau))\dot\gamma(\tau)(W(\tau)-W(t))d\tau\nonumber\\
+\nabla F(\gamma(s))(W(s)-W(t))
\label{eq:dotgamma}
\end{eqnarray}
for all $t-1\le s\le t$. Therefore
\beqas
|\dot\gamma(s)|\ge |\dot\gamma(t)|-c\max_{t-1\le\tau\le t}|W(\tau)-W(t)|\Big(1+\int_{t-1}^t|\dot\gamma(\tau)|d\tau\Big)\\
\ge |\dot\gamma(t)|-c(1+\max_{t-1\le\tau\le t}|W(\tau)-W(t)|^2).
\eeqas
From this and \eqref{eq:dotgammat}, we deduce that
\beqas
|\dot\gamma(t)|\le c(1+\max_{t-1\le\tau\le t}|W(\tau)-W(t)|^2)+c(1+\max_{t-1\le\tau\le t}|W(\tau)-W(t)|)\\
\le c(1+\max_{t-1\le\tau\le t}|W(\tau)-W(t)|^2).
\eeqas
\hfill$\Box$

We have the following lemma on the Lipschitzness of the solution $\phi^W_b(\cdot,t)$. 
\begin{lemma}\label{lem:Lipschitzp}
The solution $\phi^W_b(\cdot,t)$ satisfies
\[
|\phi^W_b(x,t)-\phi^W_b(x',t)|\le c(1+\max_{t-1\le\tau\le t}|W(\tau)-W(t)|^2)|x-x'|.
\]
\end{lemma}
The lemma is proved in \cite{IK}, but to show the explicit Lipschitz constant, we briefly mention the proof here. 

{\it Proof}\ \ Let $\gamma$ be a one-sided minimizer starting at $(x,t)$. We consider the curve
\[
\tilde\gamma(\tau)=\gamma(\tau)+(x'-x)(\tau-t+1),
\]
(modulus $\IT^d$) connecting $x'$ and $\gamma(t-1)$. As 
\[
\psi^W_b(x,t)=\psi^W_b(\gamma(t-1),t-1)+\cab_{t-1,t}(\gamma),
\]
{ and }
\[
 \psi^W_b(x',t)\le\phi^W_b(\gamma(t-1),t-1)+\cab_{t-1,t}(\tilde\gamma),
\]
we have
\beqas
&&\psi^W_b(x',t)-\psi^W_b(x,t)\le \cab_{t-1,t}(\tilde\gamma)-\cab_{t-1,t}(\gamma)\\
&&\qquad\le c\Big[1+\max_{t-1\le s\le t}|W(s)-W(t)|\Big(1+\int_{t-1}^t|\dot\gamma(\tau)|d\tau\Big)\Big]|x-x'|.
\eeqas
Similarly, the same estimate holds for $\psi^W_b(x,t)-\psi^W_b(x',t)$. 
From \eqref{eq:dotgammat}, we have
\[
|\phi^W_b(x,t)-\phi^W_b(x',t)|\le c(1+\max_{t-1\le \tau\le t}|W(\tau)-W(t)|^2)|x-x'|.
\]
\hfill$\Box$

We now establish the continuity of $\cg_B(W)$ as a map from $\cal X$ to $\IR^m$.
\begin{proposition}\label{prop:concgBp}
The function $\cg_B(W):{\cal X}\to\IR^m$ is continuous with respect to the metric $D$ in \eqref{eq:metric}.
\end{proposition}
{\it Proof}\ \ From Lemma \ref{lem:Lipschitzp}, when $W_k\to W$ in the metric $D$, $\phi^{W_k}_b$ is uniformly Lipschitz in $C^{0,1}(\IT^d)/\IR$. We can therefore extract a subsequence that converges in $C^{0,1}(\IT^d)/\IR$. From \eqref{eq:epperiodic}, every subsequence converges to $\phi^W_b$. Thus $u^{W_k}_b\to u^W_b$ whenever $\phi^{W_k}_b$ and $\phi^W_b$ are differentiable ($\phi^{W_k}_b$ and $\phi^W_b$ are differentiable almost everywhere). From Lemma \ref{lem:Lipschitzp}, $|u^{W_k}_b(x,t)|$ are uniformly bounded. The dominated convergence theorem implies that $u^{W_k}_b(\cdot,t_i)\to u^W_b(\cdot,t_i)$ in $L^1(\IT)$ for $i=1,\ldots,m$.\hfill$\Box$

For the well-posedness of the Bayesian inverse problem, we prove Assumption \ref{assumption}. 

\begin{proposition}\label{prop:assumBp}
For the spatially periodic randomly forced Burgers equation \eqref{eq:Burgers}, Assumption \ref{assumption} holds.
\end{proposition}
{\it Proof}\ \ From Lemma \ref{lem:Lipschitzp}, we deduce
\be
|\cg_B(W)|\le c\Big(1+\sum_{i=1}^m\max_{t_i-1\le\tau\le t_i}|W(\tau)-W(t_i)|^2\Big).
\label{eq:boundcgBp}
\ee
Therefore for each positive constant $M$, there exists a subset of $\cal X$ of positive Wiener measure such that for all $W$ in that set: $|\cg(W)|\le M$. By the same argument as in the proof of Proposition \ref{prop:assumHJp}, Assumption \ref{assumption}(i) holds.

From \eqref{eq:Lipy} and \eqref{eq:boundcgBp}, Assumption \ref{assumption}(ii) holds with 
\[
G(r,W)=\|\Sigma^{-1/2}\|^2_{\IR^m,\IR^m}\bigg[r+c\Big(1+\sum_{i=1}^m\max_{t_i-1\le\tau\le t_i}|W(\tau)-W(t_i)|^2\Big)\bigg],
\]
which is in $L^2({\cal X},d\mu_0)$.\hfill$\Box$

From Propositions \ref{prop:concgBp} and \ref{prop:assumBp}, Corollary \ref{corollary:2.1} and Theorem \ref{thm:wellposedness}  we have
\begin{theorem}\label{thm:Burgersp}
For the spatially periodic randomly forced Burgers equation \eqref{eq:Burgers}, the posterior measure $\mu^y$ is absolutely continuous with respect to the prior measure $\mu_0$ and satisfies 
\[
{d\mu^y\over d\mu_0}\propto \exp(-\Phi_B(W;y)),
\]
where $\Phi_B(W;y)$ is defined in \eqref{eq:PhiB}, and  is locally Lipschitz in $y$ with respect to the Hellinger metric, i.e.
\[
d_{\rm Hell}(\mu^y,\mu^{y'})\le c(r)|y-y'|,
\]
when $|y|\le r$ and $|y'|\le r$; $c(r)$ depends on $r$. 
\end{theorem}
\section{Non-periodic problems}\label{sec:nonperiodic}
We consider Bayesian inverse problems making inference on white noise forcing for  Burgers and Hamilton-Jacobi equations in non-compact domains.  Burgers and Hamilton-Jacobi equations with white noise forcing  are reasonably understood when the forcing potential $F(x)$ has a large maximum and a small minimum (Hoang and Khanin \cite{HoangKhanin}). We first present the setting up of the problem and introduce the basic results that will be used.

\subsection{Existence and uniqueness for non-periodic problems}
Let $W(t)$ be a standard Brownian motion starting at 0 at the time 0. Let 
\[
E_1=\IE\{|W(1)|\}
\]
 where $\IE$ denotes expectation with respect to the Wiener measure. For all $0<t<1$, there is a constant $E_2$ not depending on $t$ such that
\[
 \IE\{\max_{1-t\le\tau\le 1}|W(\tau)-W(1)|\}=E_2t^{1/2}.
\]
We also denote by
\[
E_3=\IE\{\max_{0\le\tau\le 1}|W(\tau)-W(1)|^2\}.
\]
Following Hoang and Khanin \cite{HoangKhanin}, we make the following assumption.
\begin{assumption}\label{assum:np}
(i) The forcing potential $F$ has a maximum at $x^{\max}$ and a minimum at $x^{\min}$ such that
\[
F(x^{\max})>L\mbox{\ \  and\ \ } F(x^{min})<-L;\ \  |x^{\max}|<a,\ \ |x^{\min}|<a
\]
 where $a$ and $L$ are positive constants; there is a constant $b>a$ such that
\[
 |\nabla F(x)|<K\ \ \mbox{when}\  |x|<b
\]
and for $x\ge b$,
\[
 |\nabla F(x)|<K_1<<K\ \mbox{ and }\ |F(x)|<L_1<<L;
\]

(ii) The function $F(x)$ is three time differentiable; $|\nabla^2F(x)|$ and $|\nabla^3F(x)|$ are uniformly bounded for all $x$. 

(iii) The constants $L,K,L_1,K_1$ satisfy 
\[
\displaystyle {128a^4E_2^2K^2\over E_1^3}<L^3,\ \ \
8a^2<LE_1,\ \ \
K_1^2<\displaystyle{LE_1\over 3E_3},\ \  
\mbox{and}\ \
L_1<\displaystyle {L\over 16}.
\]
\end{assumption}

We denote by $B(b)$ the ball with centre 0 and radius $b$. 

The solution of the Hamilton-Jacobi equation \eqref{eq:HJ}, denoted by $\phi^W$ in this case, satisfies the Lax-Oleinik formula
\be
\phi^W(x,t)=\inf_{\gamma\in\cc(s,t;\IR^d)}\{\phi(\gamma(s),s)+\ca_{s,t}(\gamma)\},
\label{eq:phinp}
\ee
where $\cc(s,t;\IR^d)$ is the space of absolutely continuous functions from $(s,t)$ to $\IR^d$; 
the operator $\ca_{s,t}$ is defined as
\begin{eqnarray}
\ca_{s,t}(\gamma)=\int_{s}^{t}\biggl({1\over 2}|{\dot\gamma}(\tau)|^2+\nabla F(\gamma(\tau)){\dot\gamma}(\tau)(W(\tau)-W(t))\biggr)d\tau+\nonumber\\
F(\gamma(s))(W(s)-W(t)).
\label{eq:AW}
\end{eqnarray}
We denote \eqref{eq:phinp} as
\[
\phi^W(x,t)={\cal K}_{s,t}^W\phi^W(\cdot,s).
\]
Minimizers for the operator $\ca_{s,t}$ are defined as in definition \ref{def:minimizer}.  

Under Assumption \ref{assum:np}, the following existence and uniqueness results for Hamilton-Jacobi equation \eqref{eq:HJ} and Burgers equation \eqref{eq:Burgers} are established in \cite{HoangKhanin}.
\begin{theorem}\label{thm:HK}
(i) With probability 1, there is a solution $\phi^W$ to the Hamilton-Jacobi equation (\ref{eq:HJ}) such that for all $t>0$, $\phi^W(\cdot,t)={\cal K}^W_{s,t}\phi^W(\cdot,s)$ for all $s<t$. The solution $\phi^W(\cdot,t)$ is locally Lipschitz for all $t$. Let $\phi_s^W(x,t)$ be the solution to the Hamilton-Jacobi equation with the zero initial condition at time $s$, then there is a constant $C_s$ such that
\[
\lim_{s\to-\infty}\sup_{|x|\le M}|\phi^W(x,t)-\phi_s^W(x,t)-C_s|=0
\]
for all $M>0$. The solution is unique in this sense up to an additive constant.
 
(ii) For every $(x,t)$, there is a one-sided minimizer $\gamma^W_{x,t}$ starting at $(x,t)$, and there is a sequence of time $s_i$ that converges to $-\infty$, and zero minimizers $\gamma_{s_i}$ on $[s_i,t]$ starting at $(x,t)$ such that $\gamma_{s_i}$ converge uniformly to $\gamma^W_{x,t}$ on any compact time intervals when $s_i\to-\infty$; further $\gamma^W_{x,t}$ is a $\phi^W$ minimiser on any finite time interval. If $x\in\spa$ is a point of differentiability of $\phi^W(\cdot,t)$ then such a one-sided minimizer starting at $x$ at time $t$ is unique and $\dot\gamma_{x,t}^W(t)=\nabla\phi^W(x,t).$

(iii) The function $u^W(\cdot,t)=\nabla\phi^W(\cdot,t)$ is a solution of the Burgers equation \eqref{eq:Burgers} that exists for all time. It is the limit of the solutions to \eqref{eq:Burgers} on finite time intervals with the zero initial condition.  In this sense, the solution $u^W$ is unique.  
\end{theorem}

From now on, by onesided minimizers, we mean those that are limits of zero minimizers on finite time intervals as in Theorem \ref{thm:HK}(ii).

The key fact that Hoang and Khanin \cite{HoangKhanin} used to prove these results is that
when $t-s$ is sufficiently large, any zero minimizers starting at $(x,t)$ over a time interval $(s,t)$ will be inside the ball $B(b)$ at a time in an interval $(\bar t(W),t)$ where $\bar t(W)$ only depends on $W$, $x$ and $t$. 
 
Let 
\[
\alpha=\min\{1,L^2E_1^2/(16a^2K^2E_2^2)\}.
\label{eq:tbar}
\] 
We define for each integer $l$ the random variable
\begin{eqnarray}
P_l={2a^2\over\alpha}+2aK\max_{l+1-\bar t\le \tau\le l+1}|W(\tau)-W(l+1)|-L|W(l)-W(l+1)|+\nonumber\\
L_1|W(l)-W(l+1)|+{K_1^2\over 2}\max_{l\le\tau\le l+1}|W(\tau)-W(l+1)|^2.\quad\quad \label{eq:Pl}
\end{eqnarray}
These random variables are independent and identically distributed. Their expectation satisfies
\be
\IE\{P_l\}<-{LE_1\over 2},
\label{eq:EPl}
\ee
which can be shown from Assumption \ref{assum:np}(iii) (see \cite{HoangKhanin} page 830).

Let $t_i'$ and $t_i''$ be the integers that satisfy  
\be
t_i\le t_i'< t_i+1\ \mbox{ and }\  t_i-2<t_i''\le t_i-1.
\label{eq:tiprime}
\ee
 The following result is shown in \cite{HoangKhanin}:
\begin{proposition}\label{prop:compact}(\cite{HoangKhanin} pages 829-831) 
Fixing an index $i$, there are constants $c(x_i,t_i,F)$ and $c(F)$ such that if ${\bar t}_i(W)$ is the largest integer for which
\begin{eqnarray}
\sum_{l=\bar t}^{t_i''-1}P_l> -c(x_i,t_i,F)\Big(1+\sum_{j=t_0''}^{t_i'}\max_{j\le\tau\le j+1}|W(\tau)-W(j+1)|^2\Big)-\nonumber\\
c(F)\Big(1+\max_{\bar t-2\le\tau\le \bar t}|W(\tau)-W(\bar t)|^2\Big),\label{eq:sumPl}
\end{eqnarray}
for all integers $\bar t\le\bar t_i(W)+2$, then when $s$ is sufficiently small, every zero minimizer over $(s,t_i)$ starting at $(x_i,t_i)$ is inside the ball $B(b)$ at a time $s\in [\bar t_i(W),t_0'']$.
\end{proposition}
Indeed, Hoang and Khanin \cite{HoangKhanin} prove the result for $s\in [\bar t_i(W),t_i-1]$ but the proof for the requirement $s\in [\bar t_i(W),t_0'']$ is identical. The law of large number guarantees the existence of such a value ${\bar t}_i(W)$. It is obvious that all onesided minimizers in Theorem \ref{thm:HK}(ii) are in $B(b)$ at a time $s\in[\bar t_i(W),t_0'']$. 

We employ these results to prove the validity of the Bayes formula and the well-posedness of the Bayesian inverse problems for Burgers and Hamilton-Jacobi equations. 
\subsection{Bayesian inverse problem for Hamilton-Jacobi equation \eqref{eq:HJ}}  
We consider Hamilton-Jacobi equation \eqref{eq:HJ} with the forcing potential $F$ satisfying Assumption \ref{assum:np}. 
To prove the validity of the Bayes formula, we show the continuity of the function 
\[
{\cal G}_{HJ}(W)=(\phi^W(x_1,t_1)-\phi^W(x_0,t_0),\ldots,\phi^W(x_m,t_m)-\phi^W(x_0,t_0)).
\]
\begin{proposition}\label{prop:concgHJnp}
The function ${\cal G}_{HJ}$ is continuous as a map from $\cal X$ to ${\mathbb R}^m$ with respect to the metric $D$ in \eqref{eq:metric}.
\end{proposition}
{\it Proof}\ \ Let $W_k$ converge to $W$ in the metric $D$ of $\cal X$. As shown in section \ref{sec:6.2}, there are constants $C_k$ which do not depend on $i$ such that
\beq
\lim_{k\to\infty}|\phi^{W_k}(x_i,t_i)-\phi^W(x_i,t_i)-C_k|=0,
\label{eq:epnonperiodic}
\eeq
for all $i=1,\ldots,m$. 
From this
\[
\lim_{k\to\infty}|(\phi^{W_k}(x_i,t_i)-\phi^{W_k}(x_0,t_0))-(\phi^{W}(x_i,t_i)-\phi^{W}(x_0,t_0))|=0.
\]
The conclusion then follows.\hfill$\Box$

To show the continuous dependence of the posterior measure on the noisy data $y$, we show  Assumption \ref{assumption}. We first prove the following lemmas.
\begin{lemma}\label{lem:1}
There is a constant $c$ that depends on $x_i,x_0,t_i,t_0$ and $F$ such that
\[
\phi^W(x_i,t_i)-\phi^W(x_0,t_0)\le c(1+\max_{t_0\le\tau\le t_i}|W(\tau)-W(t_i)|).
\]
\end{lemma}
{\it Proof}\ \ Let $\bar\gamma$ be the linear curve that connects $(x_i,t_i)$ with $(x_0,t_0)$. We have
\[
\phi^W(x_i,t_i)\le \phi^W(x_0,t_0)+{\cal A}^W_{t_0,t_i}(\bar\gamma).
\]
Since $\dot{\bar\gamma}(\tau)=(x_i-x_0)/(t_i-t_0)$, from \eqref{eq:AW}
\beqas
|{\cal A}^W_{t_0,t_i}|\le {|x_i-x_0|^2\over 2|t_i-t_0|}+\max_x|\nabla F(x)|\max_{t_0\le\tau\le t_i}|W(\tau)-W(t_i)||x_i-x_0|+\\
\max_x|F(x)||W(t_i)-W(t_0)|.
\eeqas
The conclusion then follows.\hfill$\Box$.

\begin{lemma}\label{lem:2} The following estimate holds
\[
\phi^W(x_i,t_i)-\phi^W(x_0,t_0)\ge -c(x_0,x_i,F)\sum_{l={\bar t}_i(W)}^{t_i'}(1+\max_{l\le\tau\le l+1}|W(\tau)-W(l+1)|^2),
\]
where the constant $c(x_0,x_i,F)$ depends on $x_0$, $x_i$ and $F$.
\end{lemma}
{\it Proof}\ \ We prove this lemma in \ref{sec:6.3}.

From these two lemmas, we have:
\begin{proposition}\label{prop:assHJnp} Assumption \ref{assumption} holds.
\end{proposition}
{\it Proof}\ \ From Lemmas \ref{lem:1} and \ref{lem:2},
\[
|\cg_{HJ}(W)|\le S(W),
\]
where
\[
S(W)=c+c\sum_{i=1}^m\max_{t_0\le\tau\le t_i}|W(\tau)-W(t_i)|+c\sum_{i=1}^m\sum_{l={\bar t}_i(W)}^{t_i'}(1+\max_{l\le\tau\le l+1}|W(\tau)-W(l+1)|^2),
\]
where the constants do not depend on the Brownian path $W$.

Let $I$ be a positive integer such that the subset ${\cal X}_I$  containing all $W\in{\cal X}$ that satisfy ${\bar t}_i(W)>I$ for all $i=1,\ldots,m$ has a positive probability. 
There is a constant ${\bar c}$ such that the subset containing the paths $W\in{\cal X}_I$ that satisfy 
\[
S(W)<{\bar c}
\]
has a positive probability. Assumption \ref{assumption}(i) holds.

To prove Assumption \ref{assumption}(ii), we get from \eqref{eq:Lipy}
\[
|\Phi_{HJ}(W;y)-\Phi_{HJ}(W;y')|\le {1\over 2}\|\Sigma^{-1/2}\|^2_{{\mathbb R}^m,{\mathbb R}^m}(|y|+|y'|+2|{\cal G}_{HJ}(W)|)|y-y'|.
\]
Let
\beq
G(r,W)=\|\Sigma^{-1/2}\|^2_{{\mathbb R}^m,{\mathbb R}^m}(r+S(W)).
\label{eq:nonperiodicG}
\eeq
We prove that $G(r,W)$ is in $L^2({\cal X},\mu_0)$ in \ref{sec:6.4}.\hfill$\Box$

From Propositions \ref{prop:concgHJnp}, \ref{prop:assHJnp}, Corollary \ref{corollary:2.1} and Theorem \ref{thm:wellposedness}, we get
\begin{theorem}\label{thm:HJnp}
Under Assumption \ref{assum:np}, for  Hamilton-Jacobi equation \eqref{eq:HJ} the posterior probability measure $\mu^y(dW)=\IP(dW|y)$ is absolutely continuous with respect to the Wiener measure $\mu_0(dW)$; its Radon-Nikodym derivative satisfies
\[
{d\mu^y\over d\mu_0}\propto\exp(-\Phi_{HJ}(W;y))
\]
where $\Phi_{HJ}$ is defined in \eqref{eq:PhiHJ}. For $|y|\le r$ and $|y'|\le r$, there is a constant $c(r)$ depending on $r$ such that 
\[
d_{\rm Hell}(\mu^y,\mu^{y'})\le c(r)|y-y'|.
\] 
\end{theorem}

\subsection{Bayesian inverse problem for Burgers equation}
We  consider  Burgers equation \eqref{eq:Burgers} in non-compact setting. We first prove that 
\[
\cg_B(W)=(l_1(u(\cdot,t_1)),\ldots,l_m(u(\cdot,t_m)))
\]
is continuous with respect to $W$. To this end, we prove a bound for $|\dot\gamma(t_i)|$ where $\gamma(t_i)$ is a onesided minimizer. Another bound is proved in \cite{HoangKhanin} using Gronvall's inequality. This bound involves an exponential function of the Brownian motion, which is difficult for constructing the function $G(r,W)$ in Assumption \ref{assumption}.  
\begin{lemma}\label{lem:dotgammatnp}
All one sided minimizers $\gamma$ that start at $(x,t_i)$ satisfy 
\be
|\dot\gamma(t_i)|\le c(t_i-\bar t_i(W))(1+|x|^2+\max_{\bar t_i(W)\le\tau\le t_i}|W(\tau)-W(t_i)|^2),
\label{eq:dotgammatnp}
\ee
where the constant $c$ does not depend on $x$, $t_i$ and $\bar t_i(W)$. 
\end{lemma}
{\it Proof}\ \ Let $s\in (\bar t_i(W),t_i)$ such that $\gamma(s)\in B(b)$. From  \eqref{eq:AW}, we have
\begin{eqnarray}
\ca_{s,t_i}(\gamma)\ge {1\over 2(t_i-s)}\Big(\int_s^{t_i}|\dot\gamma(\tau)|d\tau\Big)^2-c\max_{s\le\tau\le t_i}|W(\tau)-W(t_i)|\int_s^{t_i}|\dot\gamma(\tau)|d\tau\nonumber\\
-c|W(s)-W(t_i)|.\quad\quad
\label{eq:eq3}
\end{eqnarray}
We consider the linear curve $\bar\gamma$ that connects $(x,t_i)$ and $(\gamma(s),s)$. As $\dot{\bar\gamma}(\tau)\le (|x|+|b|)/(t_i-s)$,
\be
\ca_{s,t}(\bar\gamma)\le {(|x|+|b|)^2\over 2(t_i-s)}+c(|x|+b)\max_{s\le\tau\le t_i}|W(\tau)-W(t_i)|+c|W(s)-W(t_i)|.
\label{eq:eq4}
\ee
The constants $c$ in \eqref{eq:eq3} and \eqref{eq:eq4} do not depend on $x$. 
From $\ca_{s,t_i}(\gamma)\le\ca_{s,t_i}(\bar\gamma)$, solving the quadratic equation, we deduce that
\[
\int_s^{t_i}|\dot\gamma(\tau)|d\tau\le (t_i-s)c(1+|x|+\max_{s\le\tau\le t_i}|W(\tau)-W(t_i)|),
\]
 where the constant $c$ does not depend on $x$. 

The rest of the proof is similar to that of  Lemma \ref{lem:dotgammat}.\hfill$\Box$

\begin{lemma}\label{lem:Lipschitznp} For any $x$ and $x'$
\beqas
|\phi^W(x,t_i)-\phi^W(x',t_i)|\le c(t_i-\bar t_i(W))^2\Big[1+|x|^2+|x'|^2+\qquad\\\max_{\bar t_i(W)\le\tau\le t_i}|W(\tau)-W(t_i)|^2\Big]|x-x'|. 
\eeqas
The proof of this lemma is similar to that of Lemma \ref{lem:Lipschitzp}.
\end{lemma}

\begin{proposition}\label{prop:concgB}
The function $\cg_B$ is a continuous map from $\cal X$ to $\IR^m$ with respect to the metric \eqref{eq:metric}.
\end{proposition}
{\it Proof}\ \ When $W_k\to W$ in the metric \eqref{eq:metric}, $\bar t_i(W_k)$ as determined in Proposition \ref{prop:compact} converges to $\bar t_i(W)$. From lemma \ref{lem:Lipschitznp}, $\phi^{W_k}(\cdot,t_i)$ are uniformly bounded in $C^{0,1}(D)/\IR$ for any compact set $D\subset\IR^m$. Therefore, from \eqref{eq:epnonperiodic}, $\phi^{W_k}(\cdot,t_i)$ converges to $\phi^{W}(\cdot,t_i)$ in $C^{0,1}(D)/\IR$. Thus, $u^{W_k}(\cdot,t_i)$ are uniformly bounded in $D$ and converges to $u^W(\cdot,t_i)$ in $L^1(D)$. Therefore $\cg_B(W^{k})\to \cg_B(W)$.\hfill$\Box$

Now we prove the continuous dependence of $\mu^y$ on $y$. 
\begin{proposition}\label{prop:ass2Bnp}
Assumption \ref{assumption} holds.
\end{proposition}
{\it Proof}\ \ From Lemma \ref{lem:Lipschitznp}, we find that at the points where $\phi^W$ is differentiable
\beqas
&&|u(x,t_i)|\le c(t_i-\bar t_i(W))\Big[1+|x|^2+\\
&&\qquad\qquad (t_i-\bar t_i(W)+1)^2\sum_{l=\bar t_i(W)}^{t_i'-1}\max_{l\le\tau\le l+1}|W(\tau)-W(l+1)|^2\Big],
\eeqas
where $t_i'$ is defined in \eqref{eq:tiprime}.
The rest of the proof is similar to that for Proposition \ref{prop:assHJnp}.

 The function $G(r,W)$ is defined as
\begin{eqnarray}
&&G(r,W)=\|\Sigma\|_{\IR^m,\IR^m}^2\bigg[r+c\sum_{i=1}^m(t_i-\bar t_i(W))\Big(1+|x|^2+\nonumber\\
&&\qquad(t_i-\bar t_i(W)+1)^2\sum_{l=\bar t_i(W)}^{t_i'-1}\max_{l\le\tau\le l+1}|W(\tau)-W(l+1)|^2\Big)\bigg],
\label{eq:GrWB}
\end{eqnarray}
where $c$ does not depend on $x$, $t_i$ and $\bar t_i(W)$. 

The proof that the function $G(r,W)$ in \eqref{eq:GrWB} is in $L^2({\cal X},d\mu_0)$ is similar to that in section \ref{sec:6.4} for the function $G(r,W)$ defined in \eqref{eq:nonperiodicG}.
\hfill$\Box$

From Propositions \ref{prop:concgB} and \ref{prop:ass2Bnp}, Corollary \ref{corollary:2.1} and Theorem \ref{thm:wellposedness}, we deduce
\begin{theorem}\label{thm:Burgersnp}
Under Assumption \ref{assum:np}, for randomly forced Burgers equation \eqref{eq:Burgers} in non-compact setting, the posterior measure $\mu^y$ is absolutely continuous with respect to the prior measure $\mu_0$ and satisfies
\[
{d\mu^y\over d\mu_0}\propto\exp(-\Phi_B(W;y))
\]
where $\Phi_B$ is defined in $\eqref{eq:PhiB}$. For $|y|\le r$ and $|y'|\le r$, there is a constant $c(r)$ depending on $r$ such that
\[
d_{\rm Hell}(\mu^y,\mu^{y'})\le c(r)|y-y'|.
\]
\end{theorem} 

\section{Proofs of technical results}

\subsection{Proof of (\ref{eq:epperiodic})}\label{sec:6.1}

We use the concept of $\ep$-narrow places defined in \cite{IK}, in particular, the following result.
\begin{lemma}\label{lem:epnarrowperiodic}
For a sufficiently large $T$, if $s_2-s_1=T$ and 
\[
\max_{s_1\le\tau\le s_2}|W(\tau)-W(s_2)|<{1\over T^2},
\]
then for any minimizers $\gamma$ on $[s_1,s_2]$, 
\[
|\cab_{s_1,s_2}(\gamma)-{b^2T\over 2}|<\ep.
\]
\end{lemma}
For each $T>0$, for almost all $W$, ergodicity implies that there are infinitely many time intervals that satisfy the assumption of Lemma \ref{lem:epnarrowperiodic}. The proof of Lemma \ref{lem:epnarrowperiodic} can be found in Iturriaga and Khanin \cite{IK} Lemma 12. We call such an interval $[s_1,s_2]$ an {\it $\ep$-narrow place}. Indeed, the definition of $\ep$-narrow places in \cite{IK} is more general, but Lemma \ref{lem:epnarrowperiodic} is sufficient for our purpose.  

To prove (\ref{eq:epperiodic}), we need a further result from \cite{IK}.
\begin{lemma}\label{lem:psinarrow}(\cite{IK}, equation (42))
Let $[s_1,s_2]$ be an $\ep$-narrow place. Let $\gamma$ be a $\psi$ minimizer on $[s_1,t]$ where $t>s_2$ for a function $\psi$. Then for all $p\in\IT^d$,
\[
\psi(p)\ge \psi(\gamma(s_1))-2\ep.
\]
\end{lemma}
{\it Proof of (\ref{eq:epperiodic})}

Let $\gamma^{W_k}$ be a onesided minimizer starting at $(x_i,t_i)$ for the action ${\cal A}^{W_k,b}$. As $D(W_k,W)\to 0$ when $k\to\infty$, from Lemma \ref{lem:dotgammat}, $\dot\gamma^{W_k}(t_i)$ is uniformly bounded. 

Let $\gamma^W$ be a one sided minimizer for ${\cal A}^{W,b}$. 
Let $[s_1,s_2]$ be an $\ep$-narrow place of the action ${\cal A}^{W,b}$ as defined in Lemma \ref{lem:epnarrowperiodic}. When $D(W,W_k)$ is sufficiently small, $[s_1,s_2]$ is also an $\ep$-narrow place for ${\cal A}^{W_k,b}$. 

Let $\gamma_1^W$ be an arbitrary $\psi_b^W(\cdot,s_1)$ minimizer on $[s_1,t_i]$. From Theorem \ref{thm:periodic}(ii), $\gamma^W$ is a $\psi_b^W(\cdot,s_1)$ minimizer on $[s_1,t_i]$, so we have from Lemma \ref{lem:psinarrow},
\[
\psi_b^W(\gamma_1^W(s_1),s_1)\ge \psi^W_b(\gamma^W(s_1),s_1)-2\ep.
\]
On the other hand, also from Lemma \ref{lem:psinarrow}:
\[
\psi_b^W(\gamma^W(s_1),s_1)\ge \psi^W_b(\gamma_1^W(s_1),s_1)-2\ep.
\] 
Therefore
\[
|\psi_b^W(\gamma^W(s_1),s_1)- \psi_b^W(\gamma_1^W(s_1),s_1)|<2\ep.
\]
Similarly, let $\gamma_1^{W_k}$ be an arbitrary $\psi_b^{W_k}(\cdot,s_1)$ minimizer on $[s_1,t_i]$. We then have
\[
|\psi_b^{W_k}(\gamma^{W_k}(s_1),s_1)- \psi_b^{W_k}(\gamma_1^{W_k}(s_1),s_1)|<2\ep.
\] 
Fix $\gamma_1^W$ and $\gamma_1^{W_k}$. Letting $C_k=\psi_b^W(\gamma_1^{W}(s_1),s_1)-\psi_b^{W_k}(\gamma_1^{W_k}(s_1),s_1)$, we deduce
\be
|\psi_b^W(\gamma^W(s_1),s_1)- \psi_b^{W_k}(\gamma^{W_k}(s_1),s_1)-C_k|<4\ep.
\label{eq:eeq}
\ee
We have
\[
\psi_b^{W_k}(x_i,t_i)=\psi_b^{W_k}(\gamma^{W_k}(s_1),s_1)+{\cal A}_{s_1,s_2}^{W_k,b}(\gamma^{W_k})+{\cal A}_{s_2,t_i}^{W_k,b}(\gamma^{W_k}).
\]
Let $\bar\gamma$ be the curve constructed as follows. For $\tau\in [s_2,t_i]$, $\bar\gamma(\tau)=\gamma^{W_k}(\tau)$. On $[s_1,s_2]$, $\bar\gamma$ is a minimizing curve for ${\cal A}^{W,b}_{s_1,s_2}$ that connects $(\gamma^{W_k}(s_2),s_2)$ and $(\gamma^W(s_1),s_1)$. We then have
\beqas
\psi_b^W(x_i,t_i)&\le& \psi_b^W(\gamma^W(s_1),s_1)+{\cal A}_{s_1,t_i}^{W,b}(\bar\gamma)\\
&=&
\psi_b^W(\gamma^W(s_1),s_1)+{\cal A}_{s_1,s_2}^{W,b}(\bar\gamma)+{\cal A}_{s_2,t_i}^{W,b}(\gamma^{W_k}).
\eeqas
Therefore
\begin{eqnarray}
\psi_b^{W}(x_i,t_i)-\psi_b^{W_k}(x_i,t_i)\le (\psi_b^W(\gamma^W(s_1),s_1)-\psi_b^{W_k}(\gamma^{W_k}(s_1),s_1))+\qquad\nonumber\\(
\cab_{s_2,t_i}(\gamma^{W_k})-{\cal A}^{W_k,b}_{s_2,t_i}(\gamma^{W_k}))+({\cal A}_{s_1,s_2}^{W,b}(\bar\gamma)-{\cal A}_{s_1,s_2}^{W_k,b}(\gamma^{W_k})).\qquad\label{eq:eeqq}
\end{eqnarray}
When $D(W_k,W)\to 0$,  $|\dot\gamma^{W_k}(t_i)|$ are uniformly bounded so from \eqref{eq:dotgamma} we deduce that $|\dot\gamma^{W_k}(\tau)|$ and $|\gamma^{W_k}(\tau)|$ are uniformly bounded for all $s_2\le\tau\le t_i$. Therefore
\be
\lim_{k\to\infty}\cab_{s_2,t_i}(\gamma^{W_k})-{\cal A}^{W_k,b}_{s_2,t_i}(\gamma^{W_k})=0.
\label{eq:e1}
\ee
From Lemma \ref{lem:epnarrowperiodic}, as $[s_1,s_2]$ is an $\ep$-narrow place for both ${\cal A}^{W,b}$ and ${\cal A}^{W_k,b}$
\be
|{\cal A}_{s_1,s_2}^{W,b}(\bar\gamma)-{\cal A}_{s_1,s_2}^{W_k,b}(\gamma^{W_k})|<2\ep.
\label{eq:eeqqq}
\ee
From \eqref{eq:eeq}, \eqref{eq:eeqq}, \eqref{eq:e1},  \eqref{eq:eeqqq} we have 
\[
\psi_b^{W}(x_i,t_i)-\psi_b^{W_k}(x_i,t_i)-C_k\le 7\ep,
\]
when $k$ is sufficiently large. 
Similarly,
\[
\psi_b^{W_k}(x_i,t_i)-\psi_b^{W}(x_i,t_i)+C_k\le 7\ep,
\] 
when $k$ is large.
From these, we deduce 
\[
|\psi_b^{W}(x_i,t_i)-\psi_b^{W_k}(x_i,t_i)-C_k|\le 7\ep,
\]
when $k$ is sufficiently large. Hence for this particular subsequence of $\{W_k\}$ (not renumbered) 
\be
\lim_{k\to\infty}|\phi_b^{W_k}(x_i,t_i)-\phi_b^W(x_i,t_i)-C_k|=0,
\label{eq:phiconvergence}
\ee
for appropriate constants $C_k$ which do not depend on $x_i$ and $t_i$. From any subsequences of $\{W_k\}$, we can choose such a  subsequence. We then get the conclusion.  
\hfill$\Box$

\subsection{Proof of (\ref{eq:epnonperiodic})}\label{sec:6.2}
The proof uses  $\ep$ narrow places for the non-periodic setting which are defined as follows. 
\begin{definition}\label{def:epnarrownonperiodic}
An interval $[s_1,s_2]$ is called an $\ep$-narrow place for the action ${\cal A}^W$defined in \eqref{eq:AW} if there are two compact sets $S_1$ and $S_2$ such that for all minimizing curves $\gamma(\tau)$ on an interval $[s_1',s_2']\supset [s_1,s_2]$ with $|\gamma(s_1')|\le b$ and $|\gamma(s_2')|\le b$, then $\gamma(s_1)\in S_1$ and $\gamma(s_2)\in S_2$. Furthermore, if $\gamma$ is a minimizing curve over $[s_1,s_2]$ connecting a point in $S_1$ to a point in $S_2$, then $|{\cal A}_{s_1,s_2}^W(\gamma)|<\ep$.
\end{definition}
In \cite{HoangKhanin} Theorem 7, such an $\ep$ narrow place is constructed as follows. Assume that $s_2-s_1=T$ and 
\[
\max_{s_1\le\tau\le s_2}|W(\tau)-W(s_2)|<{1\over 2T}.
\]
Without loss of generality, we assume that $s_1$ and $s_2$ are integers. 
Fixing a positive constant $c$ that depends only on the function $F$ and  $b$,
the law of large numbers implies that if two integers $s_1''<s_1$ and $s_2''>s_2$ satisfy
\beqas
\sum_{l=s_2}^{s_2''-1}P_l+\sum_{l=s_1''}^{s_1-1}P_l&\ge& -c\Big[1+\max_{s_1''-2\le\tau\le s_1''}|W(\tau)-W(s_1'')|^2\Big]\\
&&\qquad\qquad -c\Big[1+\max_{s_2''\le\tau\le s_2''+2}|W(\tau)-W(s_2'')|^2\Big],
\eeqas
 then $s_1''\ge s_1^*$ and $s_2''\le s_2^*$ where $s_1^*$ and $s_2^*$ depend on the path $W$ and the constant $c$ only ($P_l$ is defined in \eqref{eq:Pl}). Let 
\[
M_1=\sup_{s_1^*-2\le\tau\le s_1}|W(\tau)-W(s_1)|,\ \ M_2=\sup_{s_2\le\tau\le s_2^*+2}|W(\tau)-W(s_2)|.
\] 
Hoang and Khanin (\cite{HoangKhanin}, pages 833-835) show that there is a continuous function $R(M_1,M_2)$ and a constant $c$ that depends on $F$ such that if
\[
{c(R(M_1,M_2)^2+1)\over T}<\ep,
\]
 then $[s_1,s_2]$ is an $\ep$-narrow place of $\ca$. By using ergodicity, they show that with probability 1, an $\ep$-narrow place constructed this way always exists and $s_2$ can be chosen arbitrarily small. When $D(W,W')$ is sufficiently small, then it is also an $\ep$-narrow place for ${\cal A}^{W'}$.

We then have the following lemma (indeed Hoang and Khanin \cite{HoangKhanin}, in the proof of their Lemma 3, show a similar result for the kicked forcing case, but the proof for white noise forcing follows the same line).
\begin{lemma}
If $\gamma$ is a $\phi(\cdot,s_1)$ minimizer on  $[s_1,t]$, such that $\gamma(s_1)\in S_1$ and $\gamma(s_2)\in S_2$ then for all $p\in S_1$
\[
\phi(p,s_1)>\phi(\gamma(s_1),s_1)-2\ep.
\]
\end{lemma}
{\it Proof of (\ref{eq:epnonperiodic})}

The proof is similar to that of (\ref{eq:epperiodic}). Assume that $[s_1,s_2]$ is an $\ep$-narrow place for both $W$ and $W_k$. Any one sided minimizers starting at $(x,t)$ where $x$ is in a compact set, when $s_2$ is sufficiently smaller than $t$, must intersect the ball $B(b)$ at a time larger than $s_2$ and at a time smaller than $s_1$ (see the proof of Proposition 2 in \cite{HoangKhanin}). Therefore $\gamma(s_1)\in S_1$ and $\gamma(s_2)\in S_2$. We then proceed as in the proof of (\ref{eq:epperiodic}) in section \ref{sec:6.1} 
\hfill$\Box$

\subsection{Proof of Lemma \ref{lem:2}}\label{sec:6.3}
In this section, we prove  Lemma \ref{lem:2}.

Consider a onesided minimizer $\gamma$  that starts at $(x,t_i)$. There is a time $s\in [{\bar t}_i(W),t_0'']$ (${\bar t}_i(W)$ is defined in Proposition \ref{prop:compact}) such that $\gamma(s)\in B(b)$. 
We then have
\[
\phi^W(x_i,t_i)=\phi^W(\gamma(s),s)+\ca_{s,t_i}(\gamma).
\]
Let $s'$ be the integer such that $1\le s'-s<2$. We then have
\be
\ca_{s,t_i}(\gamma)= \ca_{s,s'}(\gamma)+\sum_{l=s'}^{t_i''-1}\ca_{l,l+1}(\gamma)+\ca_{t_i'',t_i}(\gamma),
\label{eq:eq}
\ee
where $t_i''$ is defined in \eqref{eq:tiprime}.
From \eqref{eq:AW}, we have
\beqas
\ca_{s,s'}(\gamma)&\ge& {1\over 2|s'-s|}\Big(\int_{s}^{s'}|\dot\gamma(\tau)|d\tau\Big)^2-\\
&&\max_x|\nabla F(x)|\max_{s\le\tau\le s'}|W(\tau)-W(s')|\int_{s}^{s'}|\dot\gamma(\tau)|d\tau-\\
&&\max_x|F(x)||W(s)-W(s')|,
\eeqas
which is a quadratic form for $\int_{s}^{s'}|\dot\gamma(\tau)|d\tau$. Therefore, there is a constant $c$ that depends on $F$ such that
\[
\ca_{s,s'}(\gamma)\ge -c(1+\max_{s\le\tau\le s'}(1+|W(\tau)-W(s')|^2).
\]
Performing similarly for other terms in \eqref{eq:eq}, 
we find that there is a constant $c$ which does not depends on $x,x_0,t,t_0$ such that
\[
\ca_{s,t_i}(\gamma)\ge -c\sum_{l=s''}^{t_i'}(1+\max_{l\le\tau\le l+1}|W(\tau)-W(l+1)|^2),
\]
where $t_i'$ is defined in \eqref{eq:tiprime}, and $s''$ is the integer such that $0\le s-s''<1$.
Thus
\be
\phi^W(x_i,t_i)\ge \phi^W(\gamma(s),s)-c\sum_{l={\bar t}_i(W)}^{t_i'}(1+\max_{l\le\tau\le l+1}|W(\tau)-W(l+1)|^2).
\label{eq:ee1}
\ee
Now we consider the straight line $\gamma_1$ that connects $(x_0,t_0)$ with $(\gamma(s),s)$.  We have
\[
\phi^W(x_0,t_0)\le \phi^W(\gamma_1(s),s)+\ca_{s,t_0}(\gamma_1).
\]
Let $t_0''$ be an integer such that $1\le t_0-t_0''<2$. 
We write:
\[
\ca_{s,t_0}(\gamma_1)= \ca_{s,s'}(\gamma_1)+\sum_{l=s'}^{t_0''-1}\ca_{l,l+1}(\gamma_1)+\ca_{t_0'',t_0}(\gamma_1).
\] 
The velocity $|\dot{\gamma}_1(\tau)|\le |b-x_0|$ so there is a constant $c(x_0)$ such that 
\[
\ca_{s,t_0}(\gamma_1)\le c(x_0)\sum_{l={\bar t}_i(W)}^{t_0'}(1+\max_{l\le\tau\le l+1}|W(\tau)-W(l+1)|).
\]
Therefore
\be
\phi^W(x_0,t_0)\le \phi^W(\gamma(s),s)+c(x_0)\sum_{l={\bar t}_i(W)}^{t_0'}(1+\max_{l\le\tau\le l+1}|W(\tau)-W(l+1)|).
\label{eq:ee2}
\ee
From equations \eqref{eq:ee1} and \eqref{eq:ee2}, we have
\[
\phi^W(x_i,t_i)-\phi^W(x_0,t_0)\ge -c\sum_{l={\bar t}_i(W)}^{t_i'}(1+\max_{l\le\tau\le l+1}|W(\tau)-W(l+1)|^2).
\]
\hfill$\Box$

\subsection{Proof that the function $G(r,W)$ in \eqref{eq:nonperiodicG} is in $L^2({\cal X},\mu_0)$}\label{sec:6.4}

We now show that the function $G(r,W)$ defined in (\ref{eq:nonperiodicG}) is in $L^2({\cal X},\mu_0)$. We indeed show that
\[
\sum_{l={\bar t}_i(W)}^{t_i'}(1+\max_{l\le\tau\le l+1}|W(\tau)-W(l+1)|^2),
\]
is in $L^2({\cal X},\mu_0)$. 

Let $\ep$ be a small positive constant. Fixing an integer $k$, we consider the following events:

(A): 
\[
|\sum_{l=t_i''-m}^{t_i'}\max_{l\le\tau\le l+1}|W(\tau)-W(l+1)|^2-(t_i'-t_i''+m+1)E_3|< (t_i'-t_i''+m+1)\ep
\]
for all $m\ge k$,

(B): 
\[
|\sum_{l=t_i''-m}^{t_i''-1}\max_{l\le\tau\le l+1}|W(\tau)-W(l+1)|^2-mE_3|< m\ep
\]
for all $m\ge k$,

(C):
\[
|\sum_{l=t_i''-m}^{t_i''-1}P_l-m\IE\{P_l\}|<m\ep
\]
for all $m\ge k$. 
Assume that all of these events occur. Then from (A) and (B) for $m=k$, we get
\[
\sum_{l=t_i''}^{t_i'}\max_{l\le\tau\le l+1}|W(\tau)-W(j+1)|^2\le (t_i'-t_i''+1)E_3+(t_i'-t_i''+2k+1)\ep;
\]
using (B) for $m=k$ and $m=k+1$, we have
\[
\max_{t_i''-(k+1)\le\tau\le t_i''-k}|W(\tau)-W(t_i''-k)|^2\le (2k+1)\ep+E_3;
\]
using (B) for $m=k+1$ and $m=k+2$ we have
\[
\max_{t_i''-(k+2)\le\tau\le t_i''-(k+1)}|W(\tau)-W(t_i''-(k+1))|^2\le (2k+3)\ep+E_3;
\]
and using (C) for $m=k$ we have
\[
\sum_{l=t_i''-k}^{t_i''-1}P_l\le k(\IE\{P_l\}+\ep)\le k(-{LE_1\over 8}+\ep),
\]
(note from \eqref{eq:EPl} that $\IE\{P_l\}\le -LE_1/8$).

Now suppose that all the events (A), (B) and (C) hold for $k=t_i''-{\bar t}_i(W)$, then 
from (\ref{eq:sumPl}), for $\bar t=\bar t_i(W)$, we get
\beqas
(t_i''-{\bar t}_i(W))\Big(-{LE_1\over 8}+\ep\Big)&\ge& -c(x_i,t_i,F)\bigg\{1+(t_i'-t_i''+1)E_3\\
&&+\Big[t_i'-t_i''+2(t_i''-{\bar t}_i(W))+1\Big]\ep \bigg\}\\
&&-c(F)\bigg\{1+\Big[4(t_i''-{\bar t}_i(W))+4\Big]\ep+2E_3\bigg\}.
\eeqas
Thus
\beqas
(t_i''-{\bar t}_i(W))\bigg[-{LE_1\over 8}+\ep\Big(1+2c(x_i,t_i,F)+4c(F)\Big)\bigg]\ge\qquad\\
-c(x_i,t_i,F)\Big[1+(t_i'-t_i''+1)E_3
+(t_i'-t_i''+1)\ep \Big]\\
-c(F)(1+4\ep+2E_3).
\eeqas
When $\ep$ is sufficiently small, this only holds when $t_i''-{\bar t}_i(W)$ is smaller than an integer $N$ which does not depend on the white noise, i.e. when ${\bar t}_i(W)\ge t_i''-N$. If (\ref{eq:sumPl}) holds but ${\bar t}_i(W)<t_i''-N$, then at least one of the events (A), (B) and (C) must not hold. Now we denote by ${\cal X}_k$ the subset of $\cal X$ that contains all the Brownian paths $W$ such that at least one of the events (A), (B) and (C) does not hold. 
For each number $r$, there is a constant $C(r,\ep)$ such that
\be
\IP\{{\cal X}_k\}\le {C(r,\ep)\over k^{r}},
\label{eq:BK}
\ee
(see Baum and Katz \cite{BaumKatz} Theorem 1). Therefore
\beqas
&&\IE\Bigg\{\bigg[\sum_{l={\bar t}_i(W)}^{t_i'}\Big(1+\max_{l\le\tau\le l+1}|W(\tau)-W(l+1)|^2\Big)\bigg]^2\Bigg\}\le \\
&&\qquad\IE\Bigg\{\bigg[\sum_{l=t_i''-N}^{t_i'}\Big(1+\max_{l\le\tau\le l+1}|W(\tau)-W(l+1)|^2\Big)\bigg]^2\Bigg\}+\\
&&\qquad\sum_{k=1}^{\infty}\IE\Bigg\{\bigg[\sum_{l=t_i''-k}^{t_i'}\Big(1+\max_{l\le\tau\le l+1}|W(\tau)-W(l+1)|^2\Big)\bigg]^2{\bf 1}_{{\cal X}_k}\Bigg\}.
\eeqas
From \eqref{eq:BK}, we have
\beqas
&&\IE\Bigg\{\bigg[\sum_{l=t_i''-k}^{t_i'}\Big(1+\max_{l\le\tau\le l+1}|W(\tau)-W(l+1)|^2\Big)\bigg]^2{\bf 1}_{{\cal X}_k}\Bigg\}\\
&&\le{C(r,\ep)\over k^{r/2}}\Bigg[\IE\Bigg\{\bigg[\sum_{l=t_i''-k}^{t_i'}\Big(1+\max_{l\le\tau\le l+1}|W(\tau)-W(l+1)|^2\Big)\bigg]^4\Bigg\}\Bigg]^{1/2}\\
&&\le{C(r,\ep)\over k^{r/2}}(t_i'-t_i''+k+1)^{3/2}\Bigg[\IE\Bigg\{\sum_{l=t_i''-k}^{t_i'}\Big(1+\max_{l\le\tau\le l+1}|W(\tau)-W(l+1)|^2\Big)^4\Bigg\}\Bigg]^{1/2}\\
&&\le C(r,\ep)k^{-r/2}(t_i'-t_i''+k+1)^2.
\eeqas
When $r$ is sufficiently large, $\sum_{k=1}^{\infty}k^{-r/2}(t_i'-t_i''+k+1)^{2}$ is finite. The assertion follows. 
\hfill$\Box$
\vskip 20pt
{\bf Acknowledgement} The author thanks Andrew Stuart for suggesting the metric space framework and for helpful discussions. He also thanks Kostya Khanin for introducing him to Burgers turbulence. 
\bibliographystyle{plain}
\bibliography{inverse}

\end{document}